\documentclass[letterpaper,11pt]{amsart}
\usepackage{url}
\usepackage[pdftex]{graphicx}
\usepackage{color}
\usepackage{amssymb}
\usepackage{amsfonts}
\usepackage{enumitem}   
\usepackage[linesnumbered,lined,ruled]{algorithm2e}
\usepackage{soul}
\usepackage{subcaption}
\usepackage{todonotes}
\usepackage{booktabs}

\usepackage[pdftex]{hyperref}
\hypersetup{colorlinks=true,citecolor=black,filecolor=black,linkcolor=black,urlcolor=black}

\newtheorem{theorem}{Theorem}

\newtheorem{lemma}[theorem]{Lemma}

\newtheorem{remark}[theorem]{Remark}
\newtheorem{proposition}[theorem]{Proposition}

\newtheorem{corollary}[theorem]{Corollary}

\usepackage{dutchcal}
\newcommand{\ac}{\mathcal{a}}
\newcommand{\bc}{\mathcal{b}}
\newcommand{\cc}{\mathcal{c}}
\newcommand{\R}{\mathbb{R}}   
\newcommand{\Rb}{\ensuremath{\overline{\R}}}
\newcommand{\Rpinf}{\ensuremath{\mathbb{R}\cup\{+\infty\}}}   
\newcommand{\dom}{\operatorname{dom}}
\newcommand{\co}{\operatorname{co}}
\newcommand{\myblue}{blue}
\newcommand{\ub}{\textcolor{\myblue}{u}}
\newcommand{\vb}{\textcolor{\myblue}{v}}
\newcommand{\wb}{\textcolor{\myblue}{w}}
\newcommand{\rhob}{\textcolor{\myblue}{\rho}}
\newcommand{\alphab}{\textcolor{\myblue}{\alpha}}
\newcommand{\betab}{\textcolor{\myblue}{\beta}}
\newcommand{\gammab}{\textcolor{\myblue}{\gamma}}
\newcommand{\etab}{\textcolor{\myblue}{\eta}}
\newcommand{\chib}{\textcolor{\myblue}{\chi}}
\newcommand{\gb}{\textcolor{\myblue}{g}}

\title[Closest convex PLQ function with minimal number of pieces]{Closest univariate convex linear-quadratic function approximation with minimal number of Pieces}
\author{Namrata Kundu and Yves Lucet,\\computer science, I. K. Barber faculty of Science, UBC Okanagan}
\thanks{Corresponding author: Yves Lucet yves.lucet@ubc.ca. This work was partly supported by the second author Discovery grant RGPIN-2018-03928 (Lucet) from the Natural Sciences and Engineering Research Council of Canada (NSERC)} 
\address{Computer Science, I. K. Barber Faculty of Science
	University of British Columbia Okanagan, 3187 University Way,
	Kelowna BC V1V 1V7, Canada} 
\email{yves.lucet@ubc.ca}
\urladdr{https://cmps-people.ok.ubc.ca/ylucet/}


\date{March 31, 2024}

\begin{document}
	\begin{abstract}
		We compute the closest convex piecewise linear-quadratic (PLQ) function with minimal number of pieces to a given univariate piecewise linear-quadratic function. The Euclidean norm is used to measure the distance between functions. First, we assume that the number and positions of the breakpoints of the output function are fixed, and solve a convex optimization problem. Next, we assume the number of breakpoints is fixed, but not their position, and solve a nonconvex optimization problem to determine optimal breakpoints placement. Finally, we propose an algorithm composed of a greedy search preprocessing and a dichotomic search that solves a logarithmic number of optimization problems to obtain an approximation of any PLQ function with minimal number of pieces thereby obtaining in two steps the closest convex function with minimal number of pieces. 
  
        We illustrate our algorithms with multiple examples, compare our approach with a previous globally optimal univariate spline approximation algorithm, and apply our method to simplify vertical alignment curves in road design optimization. CPLEX, Gurobi, and BARON are used with the YALMIP library in MATLAB to effectively select the most efficient solver.
	\end{abstract}
	
	\maketitle

	\section{Introduction}\label{s:introduction}
	Finding the closest convex function is important in various disciplines where convexity is often a desirable property. In control theory, especially in model predictive control, the system's dynamics or the constraints can be nonconvex. Approximating these nonconvex parts with convex PLQ functions allows the application of convex optimization techniques, simplifying controller design \cite{RAWLINGS-09}. In the domain of signal processing and compression, approximating complex signals with simpler, piecewise representations can make them easier to analyze, transmit, and store. Reducing the number of pieces in the approximation can lead to more efficient compression algorithms. For instance, methods like wavelet compression, which are akin to PLQ representations, have seen widespread use in image compression \cite{MALLAT-09}.
	
	In machine learning, particularly in regression and classification tasks, a convex loss function can be efficiently optimized. Nonconvex loss functions can be approximated by convex PLQ functions to leverage the benefits of convex optimization \cite{SHALEV-SHWARTZ-13}. For algorithms like decision trees and certain neural networks, piecewise approximations can serve as activation functions or decision boundaries. Simplifying these approximations can accelerate training and inference, and potentially reduce overfitting \cite{GOODFELLOW-16}.
	
	In economics, utility functions, production functions, and other vital functions might not always be convex. Approximating them with convex PLQ functions can simplify the analysis and lead to more tractable models. This is especially important in game theory, where convexity plays an important role in equilibrium concepts \cite{NISAN-07}. In financial modeling, especially in risk management where risks often manifest as nonconvexities in models, finding the closest convex PLQ function can be used for portfolio optimization, where the objective is often to maximize returns while ensuring a convex risk profile \cite{MARKOWITZ-52,HULL-02}. 
	
	In areas such as quantum mechanics or thermodynamics, complex behaviors can sometimes be approximated using PLQ functions, especially when a high level of precision is not required. These approximations can simplify calculations and make certain analytical solutions possible.
 
    Beyond the above application areas, one of our motivation is to build a toolbox for manipulating convex functions and operators in computational convex analysis~\cite{GARDINER-13,GARDINER-10a,GARDINER-11a,HAQUE-18,HIRIART-URRUTY-06,LUCET-97b,LUCET-05c,LUCET-10,LUCET-06}. Due to floating point errors, composing multiple operators results in nonconvex functions in practice, even though such functions, like the convex conjugate, should always be convex in theory. This issue prevents applying any algorithm that relies on convexity. After investigating convexity tests for piecewise functions~\cite{BAUSCHKE-16,SINGH-21}, the present paper proposes optimization models with convexity constraints to compute the closest convex function.
	
	Our approach is not restricted to convexity and allows to simplify piecewise functions by approximating them with a piecewise function having a minimal number of pieces. In Section~\ref{s:numerics}, we apply our algorithm to the vertical alignment problem in road design where the computation time of calculating an optimal road design depends directly on the number of intervals (called segments in road design) composing the vertical alignment. Decreasing the number of intervals greatly improves the computation time~\cite{BEIRANVAND-17,IANNANTUONO-23,MOMO-23,AYMAN-23}.
	
\begin{figure}
		\centering
		\begin{subfigure}[t]{0.498\textwidth}
			\centering
			\includegraphics[width=.7\linewidth]{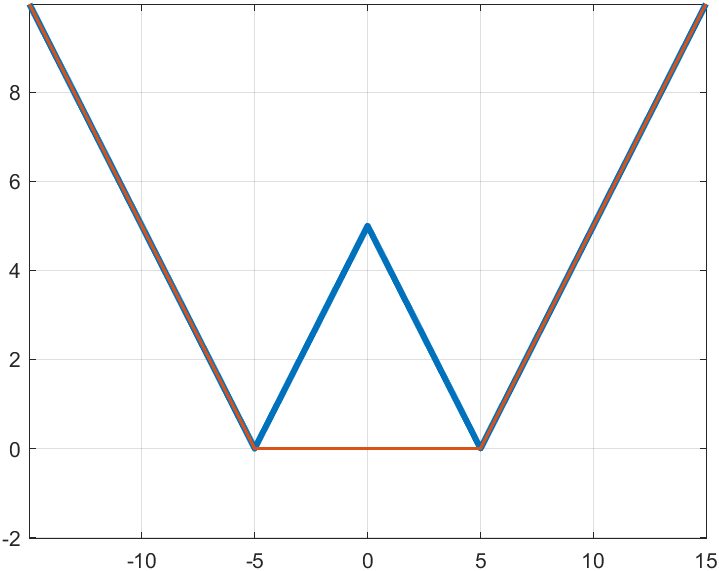}
			\caption{Convex envelope (red).}
			\label{f:convex-envelope}
		\end{subfigure}%
		~
		\begin{subfigure}[t]{0.498\textwidth}
			\centering
			\includegraphics[width=.7\linewidth]{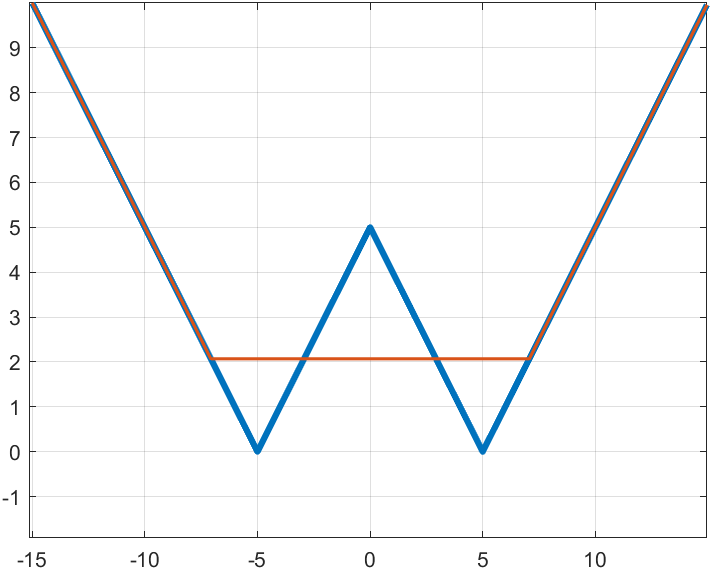}
			\caption{Closest convex function (red).}	
			\label{f:closest-convex-function}
		\end{subfigure} 
		\caption{PLQ function (blue) compared to its convex envelope and closest convex function.}	
		\label{f:covsccf}
	\end{figure}

	The closest convex function should not be confused with the convex envelope~\cite{LUCET-97b,CONTENTO-15}; see figure~\ref{f:covsccf}. We are not aware of any previous work to compute the closest convex function of a given piecewise linear-quadratic (PLQ) function. Previous work on computing piecewise functions with minimal number of intervals focused on piecewise-constant functions~\cite{KONNO-88} where the authors propose polynomial-time algorithms leveraging shortest-path and dynamic programming techniques.
	
	Piecewise-linear approximations have also been investigated with the Ramer–Douglas–Peucker algorithm~\cite{RAMER-72,DOUGLAS-73}, Visvalingam-Whyatt algorithm~\cite{VISVALINGAM-07}, and the Reumann-Witkam algorithm \cite{REUMANN-74}. The authors of~\cite{SPOERHASE-20} find that the Ramer–Douglas–Peucker runs in $O(n \log n)$ but is not optimal. For the Haussdorff distance, \cite{IMAI-88} improved by \cite{CHAN-92} is optimal in $O(n^2)$ while an optimal algorithm for the Fréchet distance runs in $O(n^3)$~\cite{GODAU-91}. Other research in this domain includes \cite{WARWICKER-21,KAZDA-24,GOLDBERG-21,REBENNACK-14,BOSCH-21}. Our numerical experiments show that our models are solved with a $O(n^2)$ complexity.
	
	Our research tackles the problem of finding the closest convex function for a given univariate PLQ function from multiple perspectives. We have formulated four optimization algorithms, each with distinct objectives and methodologies. All our algorithms are implemented in MATLAB using the YALMIP~\cite{LOEFBERG-04} numerical library so one can easily change the solver used.
	
	The first algorithm aims to determine the closest convex function with the same breakpoints as the input univariate PLQ function. The resulting convex quadratic programming optimization problem is solved effectively using CPLEX, and we obtain accurate results in polynomial time.
	
	The second algorithm deals with a more complex scenario where the number of pieces in the output convex function is set as part of the input, but their positions, that is, the breakpoints, become variable. This flexibility introduces a nonconvex nonlinear programming problem. By leveraging the state-of-the-art global optimization solver BARON~\cite{BARON-00}, we are able to attain accurate results.
	
	The above two algorithms compute the closest convex PLQ function, but may result in a representation with numerous redundant pieces. Algorithm 3 focuses on computing an optimal number of pieces to minimize the size of the function representation. It employs a dichotomy search strategy that solves an optimization problem at each step. It takes a convex PLQ function as input, which may be obtained from the output of the second algorithm. The overall objective is to determine the closest convex PLQ function with the minimum number of pieces, while adhering to a predefined error tolerance. To accomplish this, we solve an optimization problem using BARON at each step of the dichotomy search. 
	
	Lastly, we extend and generalize our algorithm 3 to approximate with minimal number of pieces any (not necessarily convex) PLQ function.
	
	Noting $n$ (resp. $m$)the number of pieces of the input (output) function, algorithm 1 (resp. 2, 3, 4) requires $n\leq m$ (resp. $n\leq m$, $n\geq m$, $n\geq m$).
	
	The paper is organized as follows. In the next section, we introduce certain definitions, notations and theorems that will be used throughout. In Section \ref{s:closest}, we formulate our first optimization models to determine the closest convex function first with fixed breakpoints and next with variable breakpoints. Section~\ref{s:minimal} focuses on minimizing the number of pieces with one algorithm preserving convexity and the other providing nonconvex approximation. Numerical experiments are reported in Section~\ref{s:numerics} including a comparison with~\cite{MOHR-23}, and an application to simplifying the representation of a vertical alignment in road design. Section~\ref{s:conclusion} summarizes our key findings and discusses avenues for future research.
	
	\section{Preliminaries and notations}\label{s:notations}
	We set our notations and recall definitions and properties used. 
	
	We note $\Rb =\R\cup\{-\infty, +\infty\}$. The (effective) domain of a function $f:\R^d\to\Rpinf$ is the set $\dom f = \{ x \in \R^d : f(x) < +\infty \}$. A function $f$ is called proper if for all $x$, $f(x)>-\infty$ and $\dom f\neq\emptyset$.
	
	A function $f: \R^d \to \Rb$ is called Piecewise Linear-Quadratic (PLQ) if $\dom f$ can be represented as the union of finitely many polyhedral sets, relative to each of which $f(x)$ is given by a quadratic expression~\cite[Definition~10.20]{ROCKAFELLAR-98a}. We call breakpoints a set $\{x_i : i=1,\dots,m\}$ with $x_i<x_{i+1}$. We allow $x_1, x_m\in \Rpinf$ while $x_i\in\R$, $i=2,\dots,m-1$. A univariate PLQ function $f$ is defined on the real line as
	\begin{equation}\label{eq:f}
		f(x) = 
		\begin{cases} 
			a_1 x^2 + b_1 x + c_1, & x \in (x_1,x_2], \\
			a_2 x^2 + b_2 x + c_2, & x \in (x_2,x_3], \\
			\vdots \\
			a_m x^2 + b_m x + c_m, & x \in (x_{m},x_{m+1}];\\
		\end{cases}
	\end{equation}
	where $a_i, b_i, c_i \in \R$ for $i=1,...,m$. To be consistent with the definition in higher dimensions, we always impose that PLQ functions are continuous on the relative interior of their domain, i.e., the only points where $f$ may be discontinuous are $x_1$ when $f(x)=+\infty$ for $x<x_1$, and $x_m$ when $f(x)=+\infty$ when $x>x_m$.
	
	We call a piece of $f$, a pair $([x_i,x_{i+1}],f_i)$ although we slightly abuse that notation to sometimes refer to an interval $[x_i,x_{i+1}]$ as a piece. 
	
	The 2-norm (also called $L_2$ norm or Euclidean norm) of a measurable function $f : \R \rightarrow \R$ is noted 
	\[
	\|f\|_2 = \sqrt{\int_{\dom f} [f(x)]^2 \, dx}.
	\]
	We say that a function $f:\R\to\Rpinf$ is coercive if $\lim_{|x| \rightarrow +\infty}f(x)=+\infty$. A function $f$ is lower semicontinuous (lsc) at a point $x$ if for every sequence $x_i \to x$, we have $f(x) \leq \liminf_i f(x_i)$.
	
	We recall the following fundamental existence result.
	
	\begin{theorem}[\cite{BAUSCHKE-24}] \label{t:key_existence_theorem}
		Let a function $f$ be a lsc coercive proper function. Then $f$ has a (global) minimizer, that is, there exists $\bar{x}$ such that
		$f(\bar{x}) = \min_{x} f(x)$.
	\end{theorem}
	
	We will also use the following corollary.
	
	\begin{corollary}[\cite{BAUSCHKE-24}] \label{c:corollary_key_existence_theorem}
		Assume $f: \R \to\Rpinf$ is coercive and lower semicontinuous. Let $C$ be a closed subset of $\R$, and assume that $C \cap \dom f \neq \emptyset$. Then $f|_C$ has a minimizer.
	\end{corollary}
	
	Finally, we call a spline a function $f:\R\to\Rpinf$ defined by a set of breakpoints $x_i$ such that (i) on each interval $[x_i,x_{i+1}]$, $f(x)$ is a polynomial of degree $k$, and (ii) $f(x)$ has continuous derivatives up to order $k-1$ at each $x_i$, for $i=2,...,m-1$.
	
	We note the following set of PLQ functions
	\[
	S(x) = \{ f : \text{$f$ is a convex PLQ function with breakpoints } x\in\R^{m}\},
	\]
	and consider the following optimization problem for a given $f\in S(x)$
	\begin{equation}\label{eq:fixed}
		\min_{\rho \in S(x)} \|f-\rho\|_2^2
	\end{equation}
	whose solution is the closest convex PLQ function that has the same breakpoints as $f$. We will also consider variable breakpoints by solving for a given $m$
	\begin{equation}\label{eq:variable}
		\min_{\eta\in\Rb^{m+1}, \rho \in S(\eta)} \|f-\rho\|_2^2.
	\end{equation}

	\section{Closest convex PLQ function}\label{s:closest}
	Given a PLQ function, we propose two algorithms to compute its closest convex PLQ function. Algorithm~1 solves~\eqref{eq:fixed}; it requires that the breakpoints of the output PLQ function $\rho$ are the same as the breakpoints of the input function $f$. By contrast, algorithm~2 only requires the number of breakpoints to solve~\eqref{eq:variable}.
	
	A function $f$ given by~\eqref{eq:f} is represented by
	\[
	\mathbf{C}(f) = \begin{bmatrix}
		a_1 & a_2 & ... & a_m \\
		b_1 & b_2 & ... & b_m \\
		c_1 & c_2 & ... & c_m \\
	\end{bmatrix}
	,
	\mathbf{P}(f) = [x_1, x_2, ... , x_{m+1}],
	\]
	where \(a = [a_1, a_2, ... , a_m]^T\), \(b = [b_1, b_2, ... , b_m]^T\), and \(c = [c_1, c_2, ... , c_m]^T\) are the vectors of the quadratic, linear, and constant coefficients of $f$ respectively; and $\mathbf{P}(f)$ is the list of breakpoints of $f$.
	
	The output of our algorithms (the closest convex PLQ function) is a PLQ function that we denote $\rho$. It is defined on $n$ pieces (algorithms 1 and 2 enforce $n \geq 2m$) and is stored as
	\begin{equation}\label{eq:rho}	
		\mathbf{C}(\rho) = \begin{bmatrix}
			\alpha_1 & \alpha_2 & ...  & \alpha_n \\
			\beta_1 & \beta_2 & ...  & \beta_n \\
			\gamma_1 & \gamma_2 & ...  & \gamma_n  \\
		\end{bmatrix}
		,
		\mathbf{P}(\rho) = [\eta_1, \eta_2,   ...  , \eta_{n+1}].
	\end{equation}
	
	In the remainder of the paper, we refer to functions defined on a piece $i$ as $f_i$, which represents $a_i x^2 + b_i x +c_i$ defined over an interval $[x_i, x_{i+1}]$.
	
	Now, we consider the boundedness of the domain of the input PLQ function $f$.
 When the input function equals a quadratic function on an unbounded interval, the closest convex function is equal to that quadratic on that interval. (If the input function has a negative quadratic coefficient on any unbounded interval, then the Euclidean distance to any convex function is infinity; it is trivial to detect and deal with this case.)
 Otherwise, we adjust the domain of the output function $\rho$ to match the domain of the input function $f$.
	
	Then Algorithm~1 returns a solution to the optimization problem (all variables are written in \textcolor{blue}{blue}).
	\begin{equation}\label{eq:algo1}
		\begin{aligned}
			\min_{\alphab_i,\betab_i,\gammab_i} \quad & \sum_{i=1}^{n} \int_{\eta_i}^{\eta_{i+1}} (f(x)- \rhob(x))^2 dx \\
			\text{subject to} \quad & \alphab_i \geq 0, 							\quad i=1,\dots,n;\\
			& \rhob_i(\eta_{i+1}) = \rhob_{i+1}(\eta_{i+1}), 		\quad i=1,\dots,n-1; \\
			& \rhob'_i(\eta_{i+1}) \leq \rhob'_{i+1}(\eta_{i+1}), \quad i=1,\dots,n-1. \\
		\end{aligned}
	\end{equation}
	All 3 constraints are linear and enforce the convexity of the output function $\rho$. The resulting optimization problem is a convex quadratic programming problem. Note that the optimization model only considers the function $\rho$ on a bounded domain making the integral in the objective function finite since all functions are continuous.
	
	Algorithm 2 sets $\rho_1$ and $\rho_n$ like algorithm 1, i.e., they are equal to the input function on any unbounded interval. It then proceeds to solve for a given tolerance $\delta>0$
	\begin{equation}\label{eq:algo2}
		\begin{aligned}
			\min_{\alphab_i,\betab_i,\gammab_i,\etab_i} \quad & \sum_{i=1}^{n} \int_{\eta_i}^{\eta_{i+1}} (f(x)-\rhob(x))^2 dx \\
			\text{subject to} \quad & \alphab_i \geq 0, 							\quad i=1,\dots,n; \\
			& \rhob_i(\eta_{i+1}) = \rhob_{i+1}(\eta_{i+1}), 		\quad i=1,\dots,n-1; \\
			& \rhob'_i(\eta_{i+1}) \leq \rhob'_{i+1}(\eta_{i+1}), \quad i=1,\dots,n-1; \\
			& \etab_i \leq \etab_{i+1} - \delta, 					\quad i=1,\dots,n-1; \\
			& \mathbf{P}(f) \subset \mathbf{P}(\rhob).
		\end{aligned}
	\end{equation}
	The constraints in \eqref{eq:algo2} include the convexity constraints similar to~\eqref{eq:algo1}, but adds that the breakpoints must be an increasing vector. It also requires the breakpoint vector $\eta$ to include the breakpoints $x$ of the input function $f$.
	
	Problem~\eqref{eq:algo2} is a nonconvex nonlinear optimization problem that requires a global solver like BARON to solve. To reduce the search space, we add the following bounds
	\begin{subequations}\label{eq:algo2constraints}
		\begin{gather}
			0 \leq \alphab_i \leq \max(a), \\
			\min(b)-N \leq \betab_i \leq \max(b)+N, \\
			\min(t_1, t_2) \leq \gammab_i \leq \max(f(x_{bounded}));
		\end{gather}
	\end{subequations}
	where $t_1$ represents the intercept of the tangent drawn at $f_1(x_{0.5})$ and $t_2$ represents the intercept of the tangent drawn at $f_m(x_{m+0.5})$ with $x_{1.5}$, $x_{m+0.5}$ as per Lemma \ref{lemma:bounding_box}. (If the conditions of the lemma are not satisfied, use $x_{1.5}=x_1$ and $x_{m+0.5}=x_{m}$.) The integer $N$ is an arbitrarily large number (default to $N=1,000$). In practice, we apply the above bounds and if we find that any bound is tight, we increase it and resolve the problem until none of the bound constraints~\eqref{eq:algo2constraints} is active.
	
\begin{figure}
		\centering
		\begin{subfigure}[t]{0.498\textwidth}
			\centering
			\includegraphics[width=.7\linewidth]{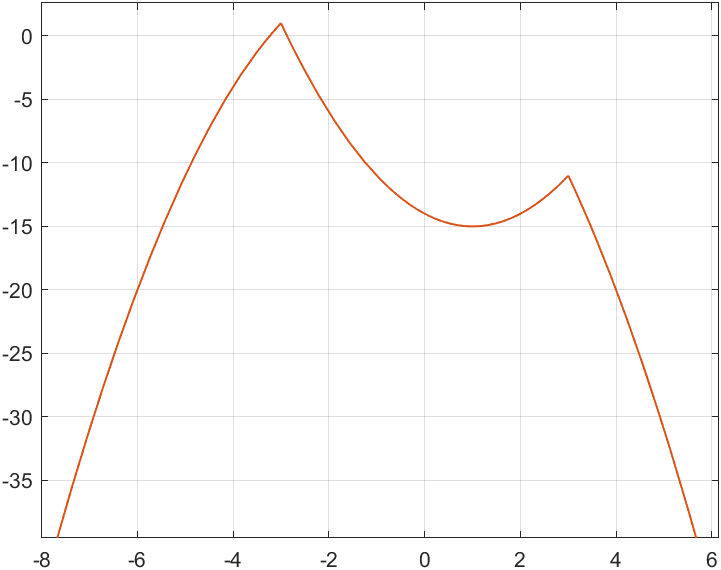}
			\caption{Case~\eqref{eq:unfeasible1}: concave quadratic.}
			\label{f:case1}
		\end{subfigure}%
		~
		\begin{subfigure}[t]{0.498\textwidth}
			\centering
			\includegraphics[width=.7\linewidth]{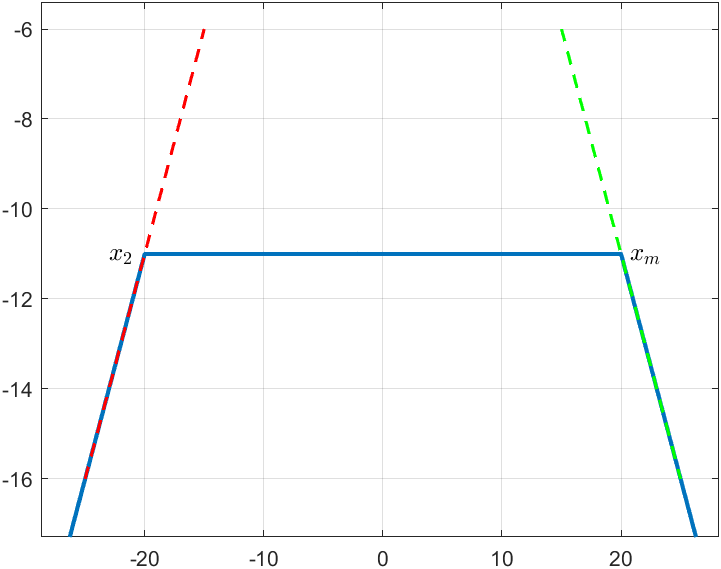}
			\caption{Case~\eqref{eq:unfeasible2}: decreasing derivative.}	
			\label{f:case2}
		\end{subfigure} 
		\caption{PLQ functions admitting no closest convex function.}
		\label{f:case12}
	\end{figure}

	We now turn our attention to the correctness of algorithms 1 and 2 by considering the following cases.
	\begin{subequations}\label{eq:unfeasible}
		\begin{gather}
			(x_1=-\infty \text{ and } a_1<0), \text{ or } (x_{m+1}=+\infty \text{ and } a_m<0),\label{eq:unfeasible1} \\ 
			x_1=+\infty, x_{m+1}=+\infty, a_1=a_m=0, \text{ and } f'_1(x_2) > f'_m(x_{m}). \label{eq:unfeasible2}		
		\end{gather}
	\end{subequations}
	In case~\eqref{eq:unfeasible1}, the input function $f$ has a strictly concave quadratic part on an unbounded interval (see figure~\ref{f:case1}), while in case~\eqref{eq:unfeasible2} the derivative values at $x_2$ and $x_{m-1}$ make it impossible to build a convex approximation whose derivative has to be nondecreasing (see figure~\ref{f:case2}).
	
	\begin{lemma}
		Assume $f$ is univariate proper PLQ. Then the following are equivalent:
		\begin{enumerate}
			\item $f$ admits a closest convex function,		
			\item neither \eqref{eq:unfeasible1} nor \eqref{eq:unfeasible2} hold,
			\item there exists an affine minorant to $f$,
			\item $\dom f^*\neq\emptyset$,
			\item $\co f \not\equiv -\infty$.
		\end{enumerate}
	\end{lemma}
	
	\begin{proof}
		If there exists an affine function minoring $f$ then $\dom f^*\neq\emptyset$; so $\co f=f^{**} \not\equiv -\infty$. The later cannot hold if either \eqref{eq:unfeasible1} or \eqref{eq:unfeasible2} hold. The last implication follow because $f$ is PLQ. Hence, (ii)--(v) hold. Lemma~\ref{lemma:case1_2} below concludes the proof.
	\end{proof}
	
	\begin{remark}
		The assumption that $f$ is PLQ is required. Consider $f(x)=1/\sqrt{-x}$ if $x\leq -1$, $f(x)=1$ if $-1<x\leq 1$ and $f(x)=1/\sqrt{x}$ if $x>1$. Then conditions (ii)--(iv) hold, but  $\|f\|=+\infty$ and $f$ does not admit a closest convex function.
	\end{remark}
	
	\begin{lemma} \label{lemma:case1_2}
		Assume f satisfies~\eqref{eq:unfeasible1} or \eqref{eq:unfeasible2}. Then $||f - g|| = \infty$ for any convex function $g$.
	\end{lemma}
	
	\begin{proof}
		For case~\eqref{eq:unfeasible1}, assume that $x_1=-\infty$ and $a_1<0$, and note $g$ any convex function. Consider the linear function $L$ that goes through $x_2$, $g(x_2)$ with slope $g'(x_2)$. Then $\infty = ||L-f|| \leq ||f-g||$. The symmetric case $x_{m+1}=\infty$ is proven similarly. 
		
		Now consider case \eqref{eq:unfeasible2}, and assume $g$ is a convex function. If $g\neq f_1$ on $(x_1,x_2]$ or $g\neq f_m$ on $[x_{m},x_{m+1})$, then $\|f-g\|=\infty$. But if $g=f_1$ on $(x_1,x_2]$ and $g= f_m$ on $[x_{m},x_{m+1})$ then $g$ cannot be convex because $g'(x_2)=f_1'(x_2) > f_m'(x_m)=g'(x_m)$.
	\end{proof}
	
	Note that in some cases, we have to modify the breakpoints to obtain a closest convex function.
	
	\begin{lemma}\label{lemma:bounding_box}
		Assume $x_1=-\infty$, $x_{m+1}=+\infty$, $f'(x_2) > f'(x_{m-1})$, and $a_1, a_{m} >0$. Then there exists $x_1 < x_{1.5} < x_2$, and $x_{m+0.5} > x_{m}$ such that $f'(x_{1.5}) \leq f'(x_{m+0.5})$.
	\end{lemma}
	
	\begin{proof}
		We have	$f'(x)\to -\infty$ when $x\to-\infty$. So there exists $x_{1.5}\in\R$ such $f'(x_{1.5}) < f'(x_2)$. Similarly, there exists $x_{m+0.5}\in\R$ with $f'(x_{m}) < f'(x_{m+0.5})$.
	\end{proof}
	
\begin{figure}
		\centering
		\begin{subfigure}[t]{0.498\textwidth}
			\centering
			\includegraphics[width=.7\linewidth]{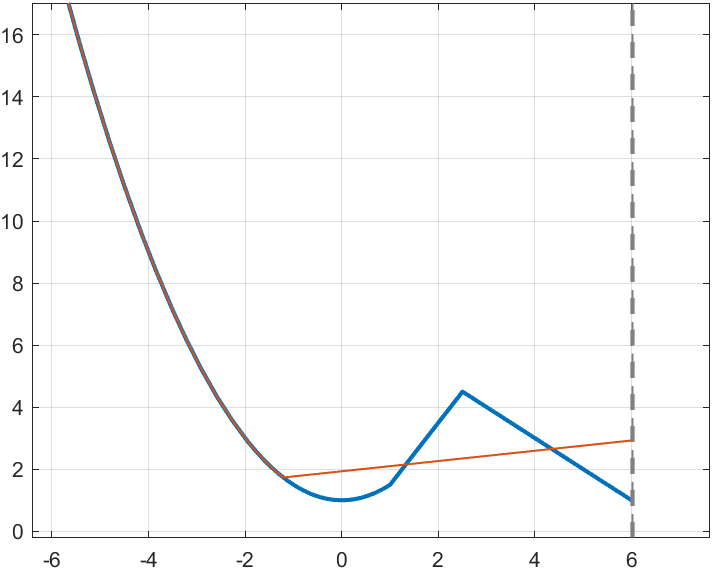}
			\caption{Added $x_{0.5}$ to breakpoint list to compute closest convex function.}
			\label{f:x0.5a}
		\end{subfigure}%
		~
		\begin{subfigure}[t]{0.498\textwidth}
			\centering
			\includegraphics[width=.7\linewidth]{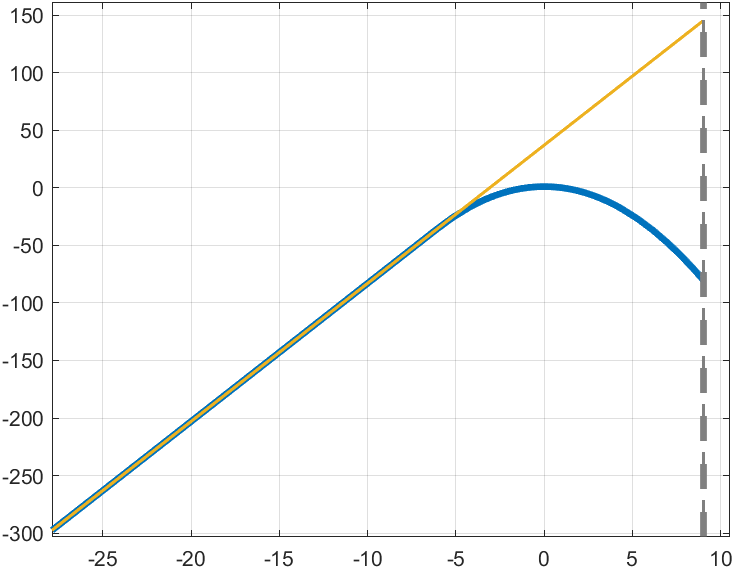}
			\caption{Concave function with right-bounded domain.}	
			\label{f:x0.5b}
		\end{subfigure} 
		\caption{PLQ functions admitting a closest convex function after applying lemma~\ref{lemma:bounding_box}.}
		\label{f:x0.5}
	\end{figure}

	Lemma \ref{lemma:bounding_box} is illustrated in figure~\ref{f:x0.5}.

	We can now replace $\mathbf{P}(f)$ with $[x_{1.5}, \mathbf{P}(f)[2:m], x_{m+0.5}]$ to obtain a different data representation for the same function $f$. With this new data representation, we can solve \eqref{eq:fixed} and \eqref{eq:variable}.
	
	\begin{proposition}\label{p:existence}
		Assume $f$ does not satisfy~\eqref{eq:unfeasible1} nor \eqref{eq:unfeasible2}. Then problems \eqref{eq:fixed} and \eqref{eq:variable} admit solutions. Moreover, algorithm~1 provides a solution to \eqref{eq:fixed}  while algorithm~2 gives a solution to~\eqref{eq:variable}.
	\end{proposition}
	
	\begin{proof}
		Invoking lemma~\ref{lemma:bounding_box} if needed, if $x_1=-\infty$ (resp. $x_{m+1}=+\infty$), algorithm 1 sets $\rho_1=f_1$ (resp. $\rho_m=f_n$) and reduces the calculation of the cost function to computing a finite integral. The resulting function is a solution to~\eqref{eq:fixed}. The fact that~\eqref{eq:fixed} has a solution is a consequence of theorem~\ref{t:key_existence_theorem}. Similarly,~\eqref{eq:variable} has (at least) one solution and algorithm~2 returns a solution. Finally, corollary~\ref{c:corollary_key_existence_theorem} shows that~\eqref{eq:variable} still admits a solution after we add constraints~\eqref{eq:algo2constraints}.
	\end{proof}
	
	\begin{remark}
		Algorithm 2 impose $n\geq 2m$ because our optimization models require $\mathbf{P}(f) \subset \mathbf{P}(\rho)$ in order to find breakpoints in $\rho$ which might fall in between existing breakpoints in $f$. Thus we change the representation of $f$ by including more breakpoints (the function $f$ remains unchanged).
	\end{remark}
	
	\begin{remark}
		For both algorithms~1 and 2, the number of constraints and decision variables increase linearly with the number of pieces $n$ in the output convex function.
	\end{remark}
	
	The optimization model formulation \eqref{eq:fixed} is a convex quadratic programming problem. We use MATLAB \cite{MATLAB-22} with YALMIP \cite{LOEFBERG-04} to interface with the CPLEX \cite{NICKEL-20} solver to compute solutions. (We also tried using Gurobi~\cite{GUROBI-24} as a solver, however, we found Gurobi was unable to solve certain cases, whereas CPLEX successfully solved for all inputs.)

\section{Minimal number of pieces}\label{s:minimal}
Algorithms 1 and 2 in the previous section require as input the number of pieces of the output function. That number is selected conservatively thereby obtaining a closest convex function representation with more pieces than necessary. This nonoptimal representation increases any computation involving such functions. Working with functions composed of fewer pieces simplifies calculations and decreases computational complexity. With this motivation, the objective of the present section is to develop an algorithm that takes as input a univariate convex PLQ function $\rho$, and outputs a univariate convex PLQ function $g$ that approximates $\rho$ within an $\varepsilon$ error tolerance and has a minimum number of piece representation. Afterwards, we drop the convexity assumption and propose an algorithm for any PLQ function.
	
Our input is stored as the output of the algorithm of the previous section, namely~\eqref{eq:rho}. Our output is a PLQ function 
 \begin{equation} \label{eq:g}
	 g(x) = \begin{cases}
		 \ac_1 x^2 + \bc_1 x +\cc_1 & \chi_0 < x \leq \chi_1, \\
		 \ac_2 ^2 + \bc_2 x +\cc_2 & \chi_1 < x \leq \chi_2, \\
		 \vdots \\
		 \ac_r x^2 + \bc_r x +\cc_r & \chi_r < x < \chi_{r+1}; \\
		 \end{cases}
	 \end{equation}
defined on $r$ pieces ($r \leq n$) and stored as
\[
\mathbf{C}(g) = \begin{bmatrix}
	\ac_1 & \ac_2  & ... & \ac_r   \\
	\bc_1 & \bc_2  & ... & \bc_r   \\
	\cc_1 & \cc_2  & ... & \cc_r   \\
\end{bmatrix}
,
\mathbf{P}(g) = [ \chi_1, \chi_{2}, ... , \chi_r, \chi_{r+1}].
\]

\begin{figure}
	\centering
	\includegraphics[width=.7\linewidth]{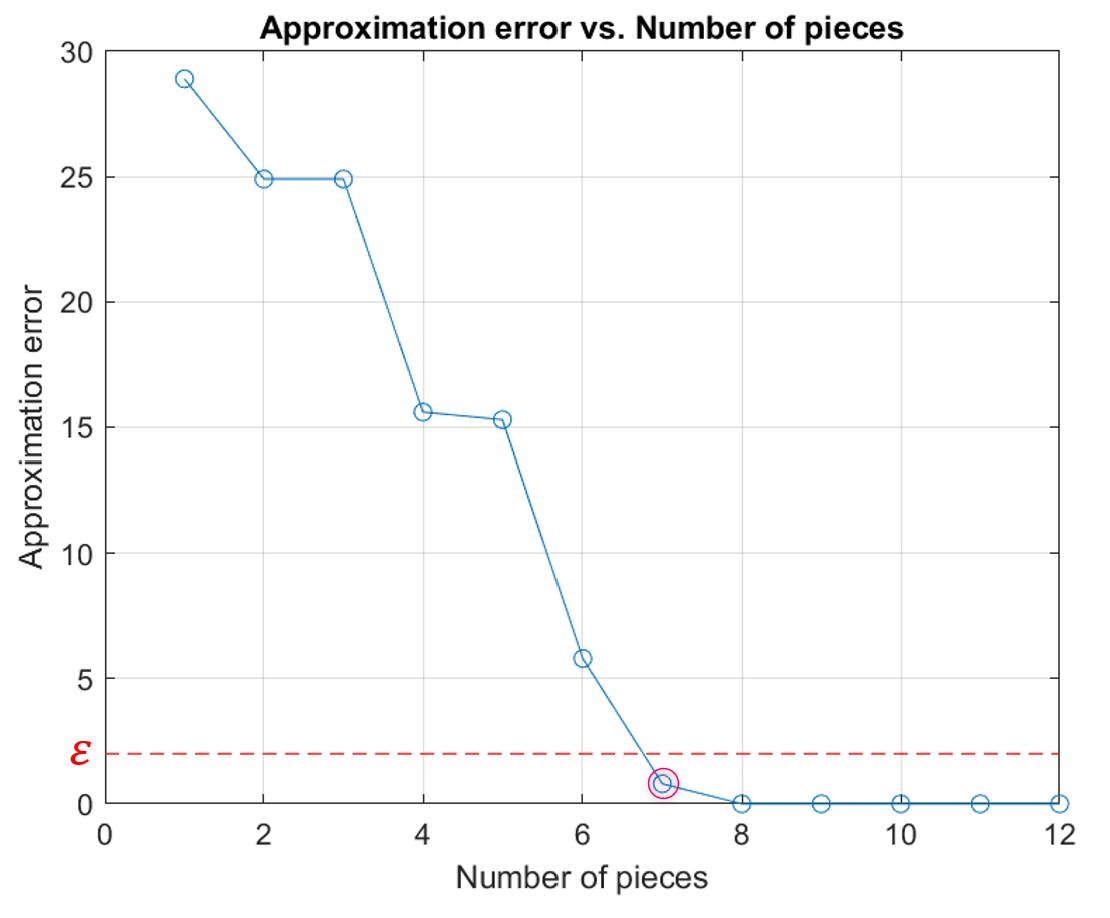}
	\caption{Approximation error vs. number of pieces. Using the nonincreasing property of the graph, we want to compute $r$ corresponding to the red point, i.e., minimal number of pieces with approximation error lower than the error tolerance $\varepsilon$.}
	\label{f:monotone-algo3}
\end{figure}

When $g$ has $r=n$ pieces, we have $g=\rho$ and the approximation error is $0$. As the number of pieces $r$ reduces, the approximation error increases. Figure~\ref{f:monotone-algo3} shows a typical example of approximation error vs. number of pieces; the key point is that the graph is nonincreasing. Our goal is to compute the smallest $r$ with approximation error below $\varepsilon$. 

Our first step is to build an algorithm that takes $r$ as input and returns a PLQ approximation with minimal error. Similar to algorithm 2, we deal with unbounded intervals first, and enforce that the output function is equal to the input function on any unbounded interval.
The same existence results apply as in the previous section, e.g., if $\alpha_n <0$ then the objective function is always equal to infinity.

Next, contrary to algorithm 2, we impose $r\leq n$ leading to the following optimization problem
	\begin{align}
		\min_{\ub,\vb,\wb,\textcolor{blue}{\ac},\textcolor{blue}{\bc},\textcolor{blue}{\cc}} \quad &  \sum_{i=1}^{n} \sum_{j=1}^{r} \ub_{i,j} \int_{\eta_i}^{\eta_{i+1}} (\rho_i - \gb_j)^2, 	\label{eq:algo3}	 \\ 
		\text{subject to} \quad &  \textcolor{\myblue}{\ac}_j \geq 0, & j=1,\dots,r; \label{eq:algo3-co1} \\ 
		& \gb_j(\chib_{j+1}) = \gb_{j+1}(\chib_{j+1}), \quad \gb'_j(\chib_{j+1}) \leq \gb'_{j+1}(\chib_{j+1}), & j=1,\dots,r; \label{eq:algo3-co2} \\ 
		& \chib_j \leq \chib_{j+1} - \delta, & j=1,\dots,r;  \label{eq:algo3-brpt1} \\
		& \mathbf{P}(\gb) \subset \mathbf{P}(\rho) &  \label{eq:algo3-brpt2} \\
		& \ub_{i,j} \in \{0,1\}, \vb_{i,j} \in \{-1,0,1\}, \wb_{i,j} \in \{0,1\} & i=1,\dots,n, j=1,\dots,r; ; \label{eq:algo3-bin}\\
		& \ub_{1,1} = \ub_{n,r} = 1, \quad \ub_{0,j}=\ub_{n+1,j}=0 &  j=1,\dots,r; \label{eq:algo3-bin-bd}\\
		& \sum_{j=1}^{r} \ub_{i,j} = 1, & i=1,\dots,n; \label{eq:algo3-bin1}\\
		& \sum_{i=1}^{n} \ub_{i,j} \geq1 & j=1,\dots,r;  \label{eq:algo3-bin2}\\ 
		& \vb_{i,j}=\ub_{i,j} - \ub_{i-1,j}, \quad \wb_{i,j}=|\vb_{i,j}|, & i=1,\dots,n+1, j=1,\dots,r; \label{eq:algo3-bin3}\\
		& \sum_i \wb_{i,j} = 2, \quad \sum_i \vb_{i,j} = 0, & j=1,\dots,r; \label{eq:algo3-bin4}\\
        & \sum_j \vb_{i,j} = 0, & i=1,\dots,n+1. \label{eq:algo3-bin5}
	\end{align}
Constraints~\eqref{eq:algo3-co1}--\eqref{eq:algo3-co2} enforce the convexity of $g$, and constraints~\eqref{eq:algo3-brpt1}--\eqref{eq:algo3-brpt2} deal with breakpoints.
Constraint~\eqref{eq:algo3-bin} defines 3 sets of binary variables with $u$ being the most important one while $v$ is a temporary variable (that the presolver will replace; it is only here for ease of exposing), and $w$ is introduced to store $w_{i,j}=|v_{i,j}|$. We want to impose a pattern on $u_{i,j}$ as shown in table~\ref{t:u} where each row sums to $1$ as enforced by~\eqref{eq:algo3-bin1}. Each column should have at least one nonzero value by~\eqref{eq:algo3-bin2} (while that may introduce redundant pieces, we can weed them out with the dichotomy search below). Moreover, each column $u_{i,j}$ has consecutive nonzero values. To enforce that constraint, we introduce $v_{i,j}=u_{i,j} - u_{i-1,j}\in\{-1,0,1\}$ as the difference between consecutive values row-wise in the $u_{i,j}$ matrix. We impose that the sums over the columns of $v$ be $0$ while sums over each column of $w_{i,j}$ add up to $2$ by~\eqref{eq:algo3-bin4}. Finally, any row $i\in\{1,\dots,n\}$ of $v_{i,j}$ has either all zeros, or a single 1 and a single -1; hence we require that their sum over the columns be zero by \eqref{eq:algo3-bin5}.

While $w_{i,j}=|v_{i,j}|$ is not a linear constraint, standard reformulation techniques already implemented in solvers BARON or Gurobi reformulate it as linear constraints by introducing new binary variables. For ease of presentation, we do not write such standard reformulation in our model.

\begin{table}
	\hfill
	\begin{subtable}[t]{0.498\textwidth}
		\centering
	\begin{tabular}{cccc}\toprule
		& $g_1$ & $g_2$ & $g_3$ \\
		\midrule
		Padding & 0 & 0 & 0\\
		$\rho_1$ & 1 & 0 & 0\\ $\rho_2$ & 1 & 0 & 0\\ $\rho_3$ & 1 & 0 & 0\\ $\rho_4$ & 1 & 0 & 0\\
		$\rho_5$ & 0 & 1 & 0\\ $\rho_6$ & 0 & 1 & 0\\ $\rho_7$ & 0 & 1 & 0\\ $\rho_8$ & 0 & 1 & 0\\
		$\rho_9$ & 0 & 0 & 1\\ $\rho_{10}$ & 0 & 0 & 1\\ $\rho_{11}$ & 0 & 0 & 1\\ $\rho_{12}$ & 0 & 0 & 1\\
		Padding & 0 & 0 & 0\\
		\bottomrule		
	\end{tabular}
	\caption{Binary variable $u_{i,j}$ including padding.}\label{t:u}
	\end{subtable}%
	\hfill
	\begin{subtable}[t]{0.498\textwidth}
		\centering
	\begin{tabular}{cccc}\toprule
	& $g_1$ & $g_2$ & $g_3$ \\
	\midrule
	& - & - & -\\ 
	$\rho_1$ & 1 & 0 & 0\\ $\rho_2$ & 0 & 0 & 0\\ $\rho_3$ & 0 & 0 & 0\\ $\rho_4$ & 0 & 0 & 0\\
	$\rho_5$ & -1 & 1 & 0\\ $\rho_6$ & 0 & 0 & 0\\ $\rho_7$ & 0 & 0 & 0\\ $\rho_8$ & 0 & 0 & 0\\
	$\rho_9$ & 0 & -1 & 1\\ $\rho_{10}$ & 0 & 0 & 0\\ $\rho_{11}$ & 0 & 0 & 0\\ $\rho_{12}$ & 0 & 0 & 0\\
	 & 0 & 0 & -1\\ 
	\bottomrule		
	\end{tabular}
	\caption{Binary variable $v_{i,j}$ associated with $u_{i,j}$.}\label{t:v}
	\end{subtable}
	\hfill
	\caption{Binary variable $u_{i,j}$ values for input function $\rho$ with $12$ pieces and output function $g$ with $3$ pieces. To correctly compute the approximation error, a pattern of contiguous ones on each column is required whereas each row should have only one nonzero value.}\label{t:int}
\end{table}

As with algorithm 2, we additionally impose bound constraints to reduce the search space
\begin{gather*}
\min(\alpha) \leq \textcolor{blue}{\ac}_j \leq \max(\alpha), \\
\min(\beta) - N \leq \textcolor{blue}{\bc}_j \leq \max(\beta)+N, \\
\min(t_1, t_2) \leq \textcolor{blue}{\cc}_j \leq \max(\rho(\eta_{bounded})). \\	
\end{gather*}

We are now in a position to give algorithm 3 to compute an approximate convex function with minimal number of pieces for a given convex PLQ function, see Algorithm 3 below. The preprocessing step on line~\ref{a:l1} merges consecutive sections into a single piece if the quadratic coefficients are close enough. In practice, it drastically cut on the number of pieces. The remaining of the algorithm is a standard dichotomy search to identify the minimal number of pieces by repeatedly solving~\eqref{algo3}. 

\SetAlgoRefName{3}	
\begin{algorithm}
	\caption{Approximate convex PLQ function with minimal number of pieces}\label{algo3}
	\SetKwData{Left}{left}\SetKwData{Right}{right}\SetKwData{Mid}{mid}\SetKwData{Error}{error}
	\SetKwFunction{Preprocess}{Preprocess} \SetKwFunction{Dichotomy}{Dichotomy} \SetKwFunction{Floor}{Floor}
	\SetKwInOut{Input}{input}\SetKwInOut{Output}{output}
	\Input{$n$, $\mathbf{C}(\rho), \mathbf{P}(\rho)$, $\varepsilon$, $\delta$}
	\Output{$r$, $\mathbf{C}(g), \mathbf{P}(g)$, \Error}
	\BlankLine
	\tcp{Merge similar pieces to reduce $n$} 
	\Preprocess($\rho$) 	\; \label{a:l1}
	\tcp{Dichotomy search}
	\Left=0; \Right=n; \Error=+$\infty$\;
	\While{\Right $>$ \Left}{
		\Mid = \Left + \Floor((\Right-\Left)/2)\;
		[C(g),P(g),error] = set $g_1$, $g_r$ and solve \eqref{algo3}\;
		\uIf{\Error $< \varepsilon$}{\Right = mid \;}
		\Else{\Left = mid\;}
	}
\end{algorithm}

\begin{remark}
	The preprocessing greedy step may give different results when traversing the input function from left-to-right and right-to-left. It might also introduce discontinuities in the resulting function. However, the optimization model within the dichotomy search maintains the continuity in the output function.
\end{remark}

\begin{remark}
	While algorithm 3 was designed to take the output of algorithm 2 as input, it can accommodate a wide range of PLQ functions, regardless of their initial convexity. In fact, convexity is not a prerequisite for the input function, only the input function restricted to unbounded intervals has to satisfy the assumptions of proposition~\ref{p:existence}. As a consequence when the closest convex PLQ function has less pieces than the input function, algorithm 3 successfully computes it. When more pieces are needed to approximate the closest convex PLQ function, one has to use first algorithm 2 and then algorithm 3.
\end{remark}

Considering our objective function \eqref{eq:algo3}, the optimization problem is mixed-integer polynomial (but not quadratic) that can be solved using BARON. Alternatively, we can introduce new variables to reformulate the problem as a quadratically-constrained quadratic program and solve it with Gurobi.

To approximate a PLQ function that is not necessarily convex, we can use the same algorithm except for removing the convexity constraint $g'_j(\chi_{j+1}) \leq g'_{j+1}(\chi_{j+1})$, and skipping any feasibility verification on unbounded intervals. The resulting problem is of a similar type, and is solved by MATLAB/YALMIP/BARON. We named this final algorithm, algorithm 4.

\begin{remark}
	The constraint~\eqref{eq:algo3-brpt2} namely $\mathbf{P}(\gb) \subset \mathbf{P}(\rho)$ may seem restrictive. If algorithm 3 or 4 is applied to the output of algorithm 2, it is not a limitation. However, in the general case, the constraints only provides an upper bound on the minimal number of pieces.
\end{remark}

\section{Numerical experiments}\label{s:numerics}
We now report numerical experiments run on algorithms 1, 2, 3, and 4. All the experiments are run on the $12^{th}$ Gen Intel(R) Core(TM) i7-12700H - 2.30 GHz, with 10 cores.

\subsection{Algorithms 1 and 2}
Consider the function
\begin{equation}\label{eq:plq}
f(x) = \begin{cases} 
	\frac{1}{2}x^2 + 1, & -\infty < x \leq 1, \\
	2x - \frac{1}{2}, & 1 < x \leq 2.5, \\
	-x + 7, & 2.5 < x \leq 6, \\
	x^2 - 5, & 6 < x < \infty.
\end{cases}
\end{equation}

\begin{figure}
	\centering
	\begin{subfigure}[t]{0.498\textwidth}
		\centering
		\includegraphics[width=.7\linewidth]{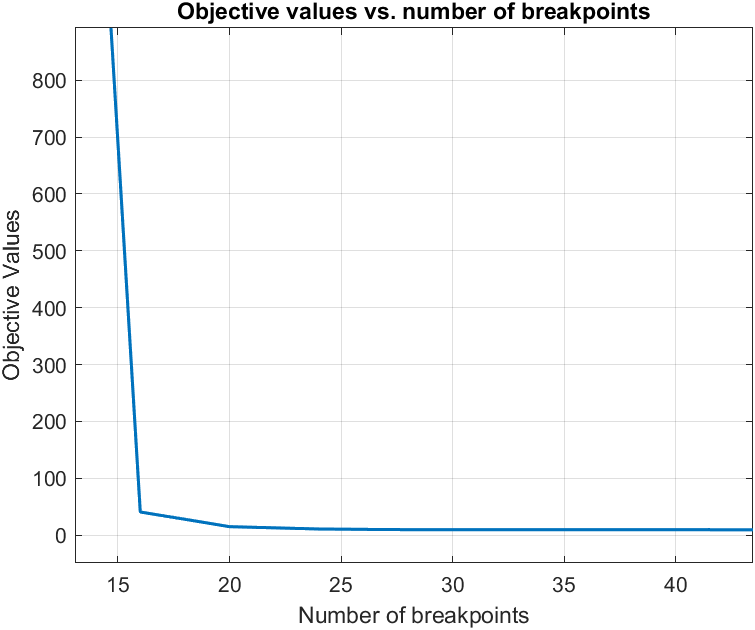}
		\caption{Distance between function and its closest convex PLQ function with the same breakpoints. Note that the objective function value is constant after 20 pieces.}
		\label{f:a1-obj-vs-bkpts}
	\end{subfigure}%
	~
	\begin{subfigure}[t]{0.498\textwidth}
		\centering
		\includegraphics[width=.9\linewidth]{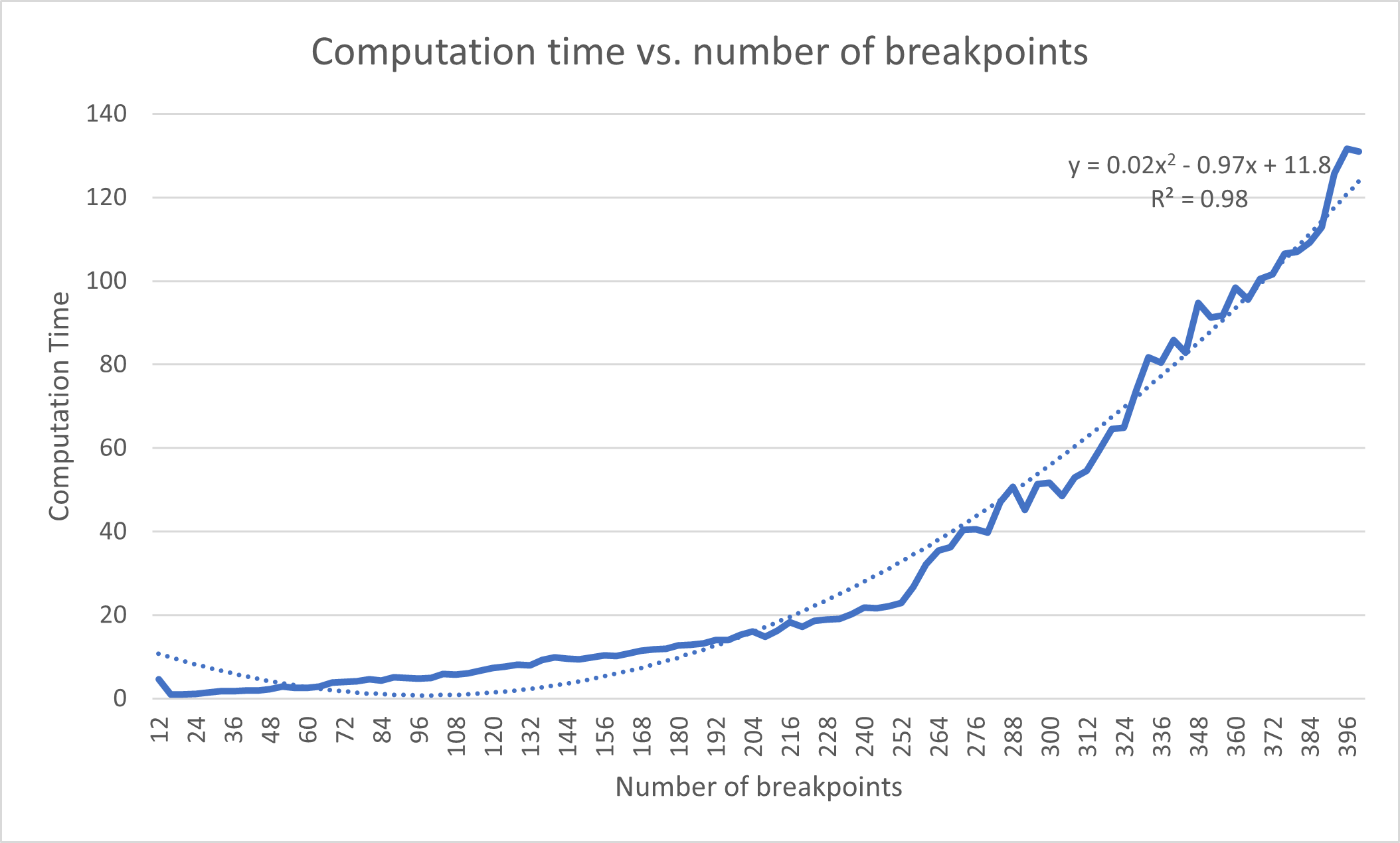}
		\caption{Computation time increases quadratically.}	
		\label{f:a1-time}
	\end{subfigure} 
	\caption{Algorithm 1.}
	\label{f:a1}
\end{figure}

We ran Algorithm 1 multiple times while increasing the number of breakpoints within each piece. As shown in figure \ref{f:a1-obj-vs-bkpts}, objective value kept decreasing till a certain number of breakpoints after which it remained constant. The computation time is plotted in figure \ref{f:a1-time}. Here computation time is the entire elapsed time to run the algorithm for a given number of output pieces - starting from problem formulation to calling the model with the solver on YALMIP. We ran the experiment thrice while clearing cache before each experiment and took the average computation time. The curve follows $O(n^2)$ increase ($n$ represents the number of breakpoints in the output function) with $R^2 = 0.98$.

\begin{figure}
	\centering
	\includegraphics[width=.7\linewidth]{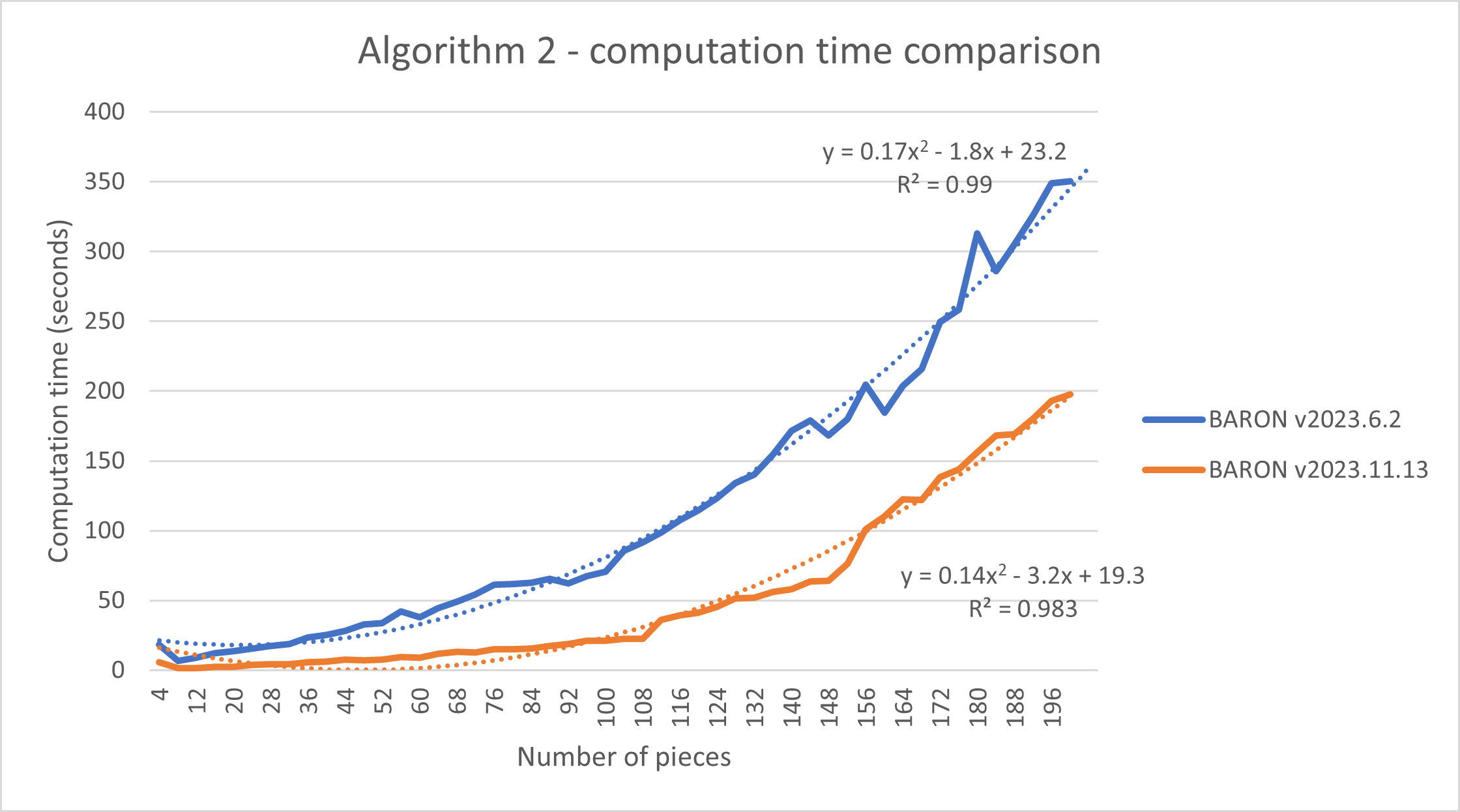}
	\caption{Algorithm 2 computation time increases quadratically. While a more recent BARON version performed better, it did not change the quadratic trend.}
	\label{f:a2-time}
\end{figure}

We did similar profiling for algorithm 2 with the same input $f$; the computation time is shown in figure~\ref{f:a2-time}. Again, we found the complexity of this algorithm grows in $O(n^2)$ with $R^2 = 0.99$. We include timings for 2 versions of BARON to emphasize that the quadratic behavior remains, but the coefficient of the quadratic curve improved with the more recent versions.

Considering that algorithms 3 and 4 are very similar, we only report numerical tests on algorithm 4.

\subsection{Comparing with globally univariate spline approximation}
We slightly modify the MIQCP Spline method from~\cite{MOHR-23} and compare it with a slight modification of algorithm 4. We adapted the MIQCP Spline method to yield piecewise quadratic splines instead of its original piecewise cubic splines, and we added a differentiability constraint at the breakpoints to algorithm~4 to ensure its output is a quadratic spline instead of a continuous piecewise quadratic polynomial, i.e., the output is $C^1$.

It is worth noting the differences between the 2 algorithms:
\begin{enumerate}
	\item algorithm 4 requires $\textcolor{blue}{\mathbf{P}(g)}\subset \mathbf{P}(\rho)$ to confine the breakpoints of the output function to be a subset of the input function's breakpoints; this constraint is not present in the MIQCP spline algorithm.
	\item the MIQCP Spline algorithm is implemented in Python and calls Gurobi while algorithm 4 is implemented in MATLAB using YALMIP and calling BARON.
	\item the MIQCP Spline algorithm minimizes the least-squares error between the input data points and the resultant function. On the other hand, algorithm~4 minimizes the area difference between the input and the output function using the 2-norm.
\end{enumerate}

Both algorithms were executed on two distinct functions: a piecewise linear function 
\[
g(x) = 
\begin{cases} 
    -x - 5, & -22 < x \leq -7.071, \\
    2.071, &-7.071 < x \leq 0.071, \\
    x - 5, &  0.071 < x < 22,
\end{cases}
\]
and the closest convex function to the PLQ function \eqref{eq:plq}.

To generate data points for the MIQCP spline algorithm, we sampled 3 data points per piece. Both algorithms accept an identical number of predetermined breakpoints (called ``knots'' in~\cite{MOHR-23}) and return a piecewise quadratic spline.

\begin{figure}
	\centering
	\begin{subfigure}[t]{0.498\textwidth}
		\centering
		\includegraphics[width=.7\linewidth]{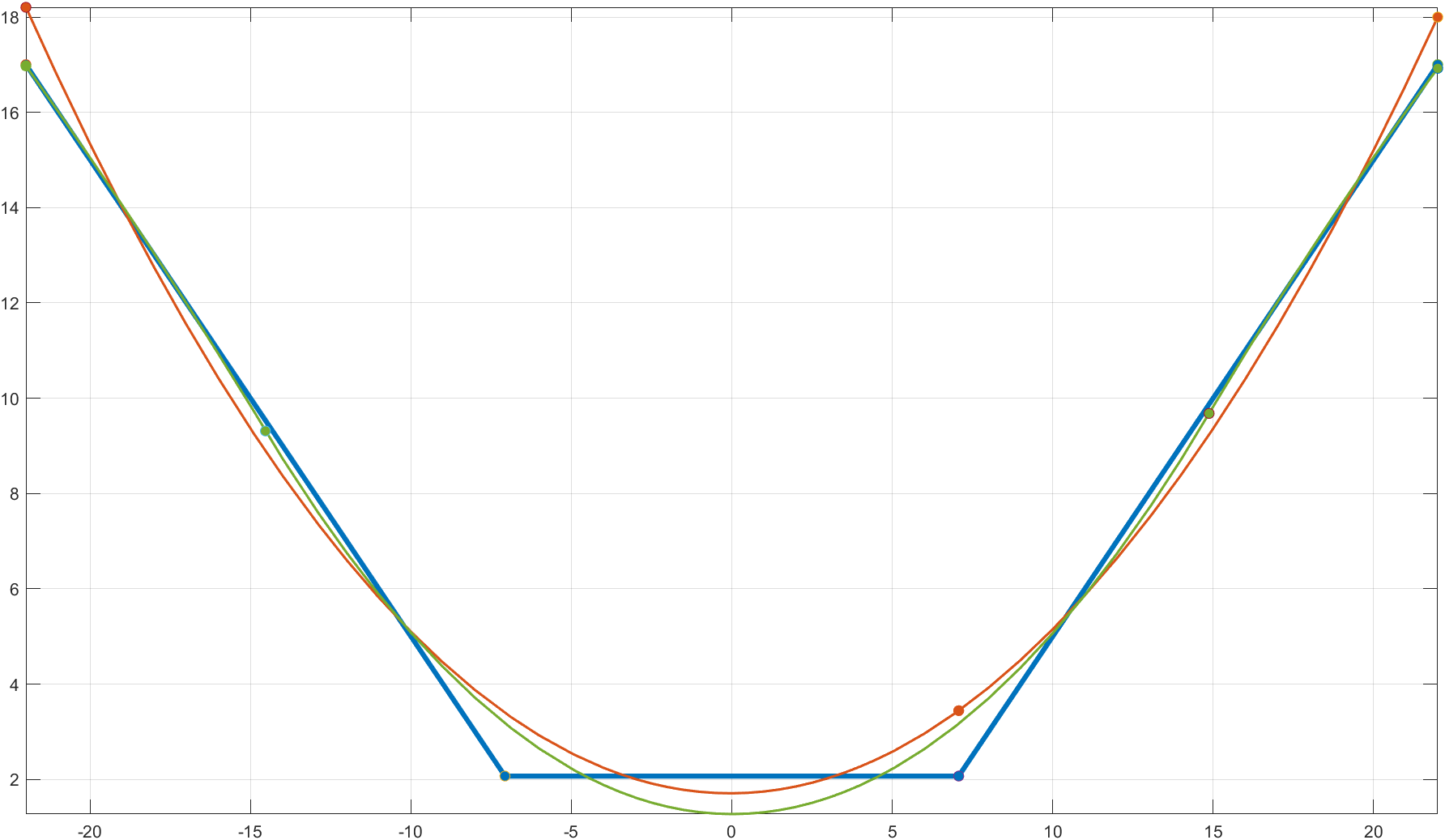}
		\caption{Approximating a piecewise linear function.}
		\label{f:a4-PL}
	\end{subfigure}%
	~
	\begin{subfigure}[t]{0.498\textwidth}
		\centering
		\includegraphics[width=.7\linewidth]{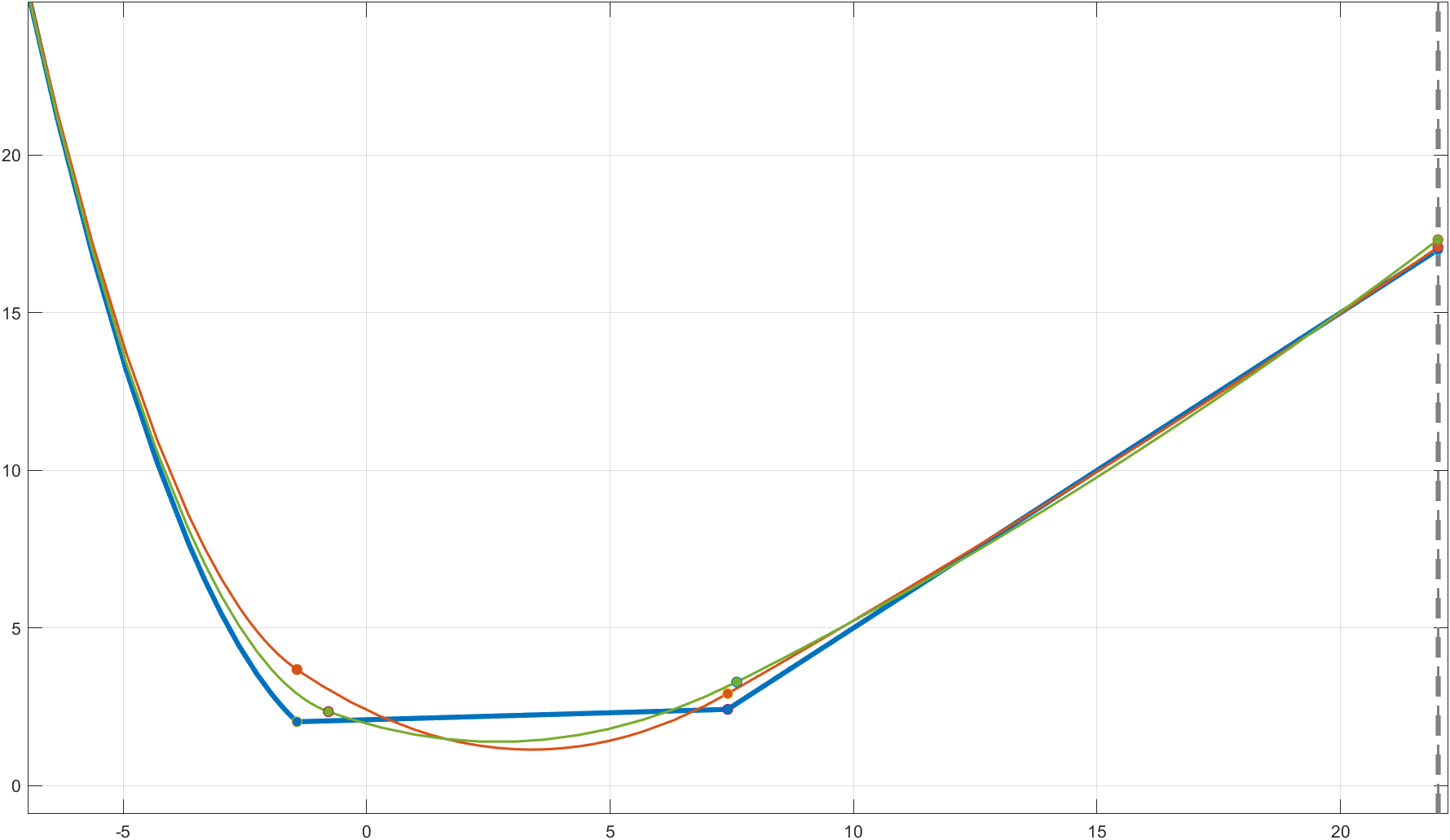}
		\caption{Approximating a piecewise quadratic function.}	
		\label{f:a4-PLQ}
	\end{subfigure} 
	\caption{Algorithm 4. Input function is shown in blue, algorithm 4 quadratic spline approximation in red, and MIQCP spline algorithm in green.}
	\label{f:a4-PL-PLQ}
\end{figure}

Figure \ref{f:a4-PL} plots the quadratic spline approximations produced by the respective algorithms. The input piecewise linear function $g$ in blue was partitioned into $36$ segments. Algorithm 4 produces a quadratic spline with $3$ pieces and a 2-norm error of $11.6$ with a computation time of $0.6$ seconds. For the MIQCP Spline algorithm, a representation of the function was created by generating $3$ data points for each piece, culminating in a set of $71$ synthetic data points. This data is combined with a specification of $2$ knots to obtain the quadratic spline in green, which has a 2-norm error of $7.7$ and requires $4.2$ seconds of computation.

We also conducted experiments using the PLQ input function $f$, with the resulting quadratic spline approximations displayed in figure \ref{f:a4-PLQ}. The input PLQ function in blue was segmented into $12$ pieces. Providing this input to the adapted algorithm 4 yielded the quadratic spline depicted in red, accompanied by a 2-norm error of $12.9$ and a computational duration of $2.7$ seconds. To adapt the function for the MIQCP Spline algorithm, $23$ synthetic data points were generated, with each piece represented by $3$ data points. Utilizing this set of data points, along with a stipulation of $2$ knots, the MIQCP Spline algorithm produced the quadratic spline represented in green, characterized by a 2-norm error of $6.1$. The computation for this problem took $0.3$ seconds.

The two algorithms return different output since they optimize different objective functions. In addition, the constraint $\textcolor{blue}{\mathbf{P}(g)}\subset \mathbf{P}(\rho)$ makes algorithm 4 suboptimal.

\begin{remark}
	As we experimented with increasing the number of pieces or breakpoints in the input function, we observed differences between solvers and between different versions of the same solver. A previous version of BARON prematurely terminated, deeming the model infeasible while the open-source solver BMIBNB for the identical model produces a solution albeit with a significantly longer computational duration (measured in hours) compared to BARON. The most recent version of BARON deemed the problem feasible. We also found discrepancies between a quadratically-constrained reformulation solved by Gurobi and BARON. Our experience suggests to use the very latest version of solvers, expect significant improvements between solver versions, and carefully validate any output.
\end{remark}

The MIQCP Spline algorithm excels in spline approximations, while algorithm 4 is specifically designed for optimizing PLQ function. Each algorithm is tailored for distinct objectives.

\subsection{Application to road design}
Designing a road involves solving three inter-related problems: drawing the horizontal alignment (a satellite view of the road), the vertical alignment (a quadratic spline drawn on a vertical plane whose $x$-axis is the distance from the start of the road, and the $y$-axis is the altitude), and the earthwork problem (indicate what material to move where to transform the ground into the vertical alignment). Our main interest is that the vertical alignment is a quadratic spline whose computation may result in a large number of pieces. This representation has an acute impact in the computation time of the road design. Consider that a road is divided into sections of length $10$ to $100$ meters. Each section may give rise to a piece in the vertical alignment and roads with hundreds of sections are common while roads of thousands of sections are becoming less rare due to the increase in computer performance. In addition, all computation depends on the precision of the ground measurements, which includes not only altitude but also estimations of various materials (rock, sand, dirts, etc.) volumes under the ground. Typically, optimization accounts for a $5$\% measurement error so there is no point in optimizing to $10$ decimals.

Algorithm 4 is well suited to approximate vertical alignment curves. We only need to add the constraint
\[\textcolor{blue}{h}'_j(\chib_{j+1}) = \textcolor{blue}{h}'_{j+1}(\chib_{j+1})\]
to ensure that the road profile is smooth. This linear constraint has no impact on the nature of the optimization problem and little impact on the computation time. 

\begin{figure}
	\centering
	\includegraphics[width=.7\linewidth]{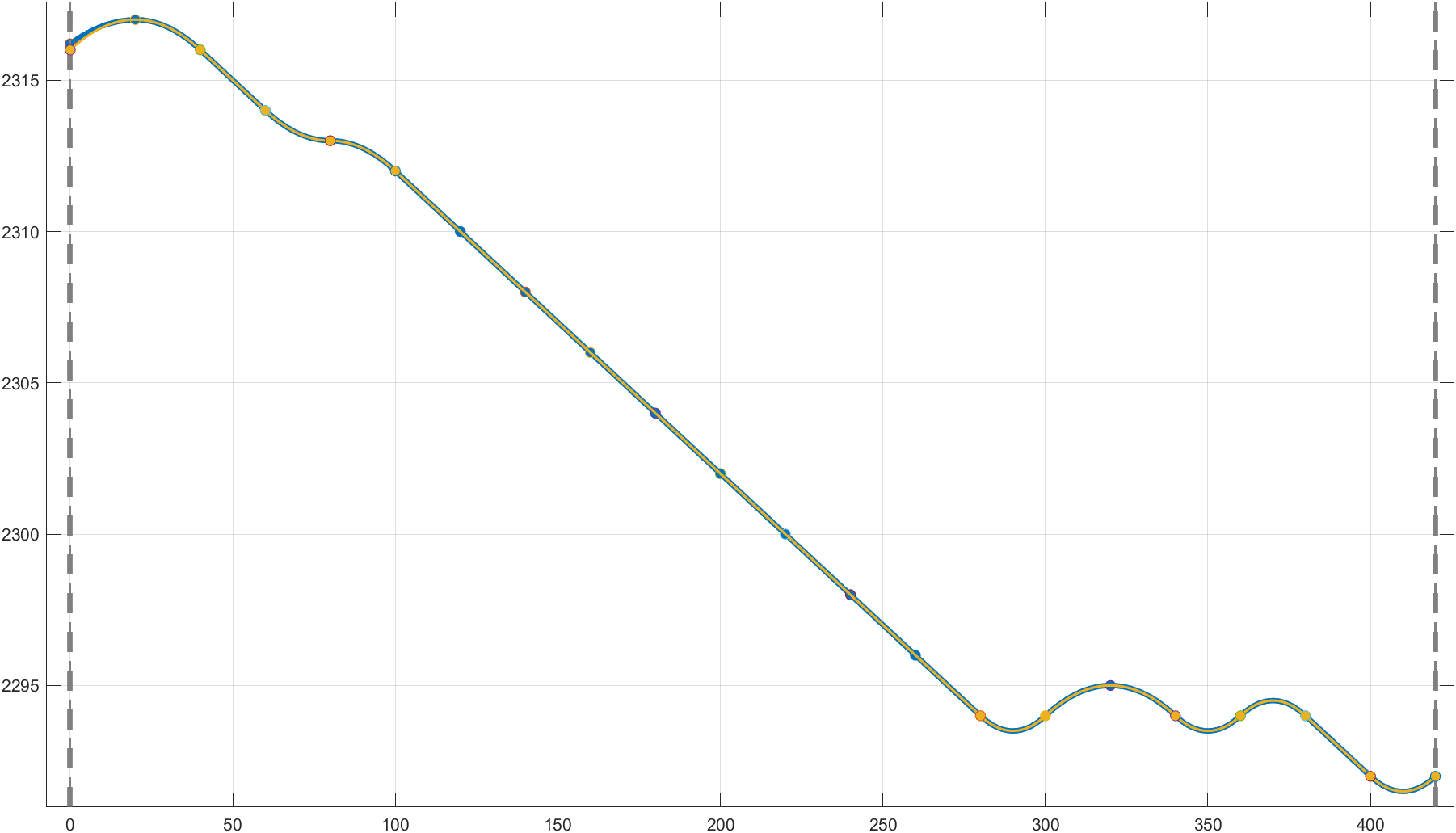}
	\caption{The input road spline $f$ is in blue and the output approximate spline function $h$ with minimal number of breakpoints in yellow. The breakpoints on each function are shown with circles in the same color. $f$ has 22 breakpoints and $h$ has 12 breakpoints.}
	\label{f:a4-VA}
\end{figure}

The data for road design was provided by Softree Technical Systems Inc.\footnote{\url{http://www.softree.com}}  
We apply algorithm~4 along with the new differentiability constraint to the vertical alignment road spline data to approximate it with the minimal number of breakpoints with an $\varepsilon$ error tolerance. The result are visualized in Figure \ref{f:a4-VA}. The input quadratic spline $f$ has $22$ breakpoints compared to $12$ breakpoints in the output quadratic spline $h$. We use $\varepsilon=0.2$ and $\delta=0.000005$; the output spline is approximated at a distance of $0.1520$ from the input spline with a computation time of $35.1$ seconds (average of 3 runs) using the quadratically-constrained quadratic programming reformulation with the Gurobi solver. 

\section{Conclusion}\label{s:conclusion}
We propose 4 algorithms. Two algorithms compute the closest convex PLQ function, and two approximate a PLQ function with a minimal number of pieces. We model optimization problems in MATLAB using YALMIP and use CPLEX, Gurobi or BARON to solve them. Numerical experiments validate the models, compare its performance to an optimal spline algorithm, and detail an application to road design optimization that greatly benefit from the reduction in the number of pieces.

More work is needed to understand differences in output between the various solvers. Some solvers sometime find a problem infeasible while other solvers can solve it. Resolving this discrepancy could expedite algorithm execution.

Since the optimization problems are polynomial, it would be interesting to substitute BARON with Gloptipoly \cite{HENRION-03}, and measure any performance difference. 

Preliminary results in reformulating our polynomial problem into a quadratically-constrained problem and solving it with Gurobi show potential for significant time reductions in addition to using a solver with a free academic license. More work is needed to systematically validate the results.

Our algorithms rely heavily on solving optimization models. Future research avenues might explore devising the closest convex function via heuristic strategies such as greedy algorithms or dynamic programming techniques. Future research could also explore parallelization techniques, e.g., to speed up the dichotomy search.

All the algorithms are specifically tailored for univariate PLQ functions. A logical progression would involve extending these algorithms to bivariate PLQ functions.

The complete code for all the algorithms is available on Github\footnote{\url{https://github.com/namrata-kundu/closest-convex-function-approximation}} under the GPL License.

\bibliographystyle{siamurl} 
\bibliography{main}

\end{document}